\documentclass[12pt]{article}

      \usepackage{latexsym}
         \usepackage[reqno, namelimits, sumlimits]{amsmath}
         \usepackage{amssymb, amsfonts}
         \usepackage{amsthm}

 \newtheorem{theorem}{Theorem}[section]

 \newtheorem{lemma}[theorem]{Lemma}
 
 \newtheorem{remark}[theorem]{Remark}

 \newtheorem{pro}[theorem]{Proposition}

\title{ Liouville Theorem for 2D Navier-Stokes equations in Half Space
}

 \author{G. Seregin
 }

\begin{document}
\maketitle

\setcounter{equation}{0}
\section{Motivation}
In the paper, we deal  with the so-called mild bounded ancient solutions to the 2D Navier-Stokes equations in half-space with the homogeneous Dirichlet boundary conditions. As it has been explained in \cite{KNSS}, \cite{SS}, and \cite{SS1}, such type of solutions appears as a result of re-scaling solutions to the Navier-Stokes equations around a possible singular point. If they are in a sense ``trivial", then this   point is not singular.

There are several interesting cases for which Liouville type theorems for ancient solutions to the Navier-Stokes equations turn out to be true. And their proofs are based on a reduction to a scalar equation with the further application  of the strong maximum principle to it. For example, in 2D case, such a scalar equations is just the 2D vorticity equation. Unfortunately, this approach does not work in a half plane since non-slip boundary conditions in terms of the velocity does not implies the homogeneous Dirchlet boundary condition for the vorticity. However, there are some interesting results coming out from this approach, see paper \cite{Giga}  and reference in it.

In the paper, we exploit a different approach related to the long time behaviour of solutions to a conjugate system. It has been already used in the proof of the Liouville type theorem for the Stokes system in half-space, see the paper \cite{JSS2} and the paper \cite{JSS1} for another approach.

Let $u$ be a mild bounded ancient solution to the Navier-Stokes equations in a half space, i.e., $u\in L_\infty(Q_-^+)$ ($|u|\leq 1$ a.e. in $Q_-^+=\{x\in\mathbb R^2_+, \,t<0\}$, where $\mathbb R^2_+=\{x=(x_1,x_2)\in\mathbb R^2:\,\,x_2>0\}$) and there exists a scalar function $p$ such that, for any $t<0$, $p=p^1+p^2$, where
\begin{equation}\label{pressure1equation}
  \triangle p^1=-{\rm div}\,{\rm div} \,u\otimes u
\end{equation}
in $Q_-^+$ with $p^1_{,2}=0$ and $p^2(\cdot,t)$ is a harmonic function in $\mathbb R^2_+$ whose gradient obeys the inequality
\begin{equation}\label{p2bounded}
  |\nabla p^2(x,t)|\leq c\ln (2+1/x_2)
\end{equation}
for   $(x,t)\in Q_-^+$ and has the property
\begin{equation}\label{pressure2property}
  \sup\limits_{x_1\in\mathbb R}|\nabla p^2(x,t)|\to 0
\end{equation}
 as $x_2\to\infty$; $u$ and $p$ satisfy the classical Navier-Stokes system and boundary condition $u(x_1,0,t)=0$
 in the following weak sense
\begin{equation}\label{weakformNSS}
 \int\limits_{ Q_-^+}\Big(u\cdot(\partial_t\varphi+\triangle\varphi)+u\otimes u:\nabla\varphi+p\,{\rm div}\,\varphi\Big)dx dt=0
\end{equation}
for any  $\varphi\in C^\infty_0(Q_-)$ with $\varphi(x_1,0,t)=0$ for $x_1\in \mathbb R$ and
\begin{equation}\label{weakdivergencefree}
   \int\limits_{ Q_-^+}u\cdot \nabla q dx dt=0
\end{equation}
for any $q\in C^\infty_0(Q_-)$.

Here, $Q_-=\mathbb R^2\times\{t<0\}$.

We are going to prove the following fact:
\begin{theorem}\label{2dLiouville} Let $u$ be a mild bounded ancient solution to the Navier-Stokes equations in a half space. Assume in addition that
\begin{equation}\label{kineticenergy}
  u\in L_{2,\infty}(Q_-^+).
\end{equation}
Then $u$ is identically equal to zero.
\end{theorem}
\begin{remark}\label{remarkaboutadditional} Motivation for additional condition (\ref{kineticenergy}) is as follows. The norm of the space $L_{2,\infty}(Q_-^+)$  is invariant with respect to the Navier-Stokes scaling
$$v(x,t)\to \lambda u(\lambda x, \lambda^2t).$$ So, if we study the smoothness of energy solutions in 2D, the corresponding norm stays bounded under scaling and a limiting procedure, leading to a mild bounded ancient solution, and thus condition (\ref{kineticenergy}) holds. For details, see \cite{SS1}. \end{remark}
\begin{lemma}\label{dissipationisbounded} Under assumptions of Theorem \ref{2dLiouville}, 
\begin{equation}\label{dissipation}
  \nabla u\in L_2(Q^+_-).
\end{equation}
 \end{lemma}
\textsc{Proof}
For fixed $A<0$, we can construct $\widetilde{u}$ as a solution
to the initial boundary value problem
$$\partial_t\widetilde{u}-\triangle\widetilde{u}+\nabla\widetilde{p}^2= -{\rm div H}$$
in $\mathbb R^2_+\times ]A,0[$, where  $H=u\otimes u+p^1\mathbb I$,
$$\widetilde{u}(x_1,0,t)=0,$$
$$\widetilde{u}(x,A)=u(x,A)$$
with the help of the Green function $G$ and the kernel $K$ introduced by Solonnikov in \cite{Sol2003}, i.e.,
$$\widetilde{u}(x,t)=\int\limits_{\mathbb R^2_+}G(x,y,t-A)u(y,A)dy+\int\limits^t_A\int\limits_{\mathbb R^2_+}K(x,y,t-\tau)F(y,\tau)dy d\tau.$$
For the further details, we refer the reader to the paper \cite{SS1}.

Let us describe the properties of $\widetilde{u}$. Our first observation is that
$$
{\rm div}\,u\otimes u=u\cdot \nabla u\in L_{2,\infty}(Q^+_-)
$$
since $u\in L_{2,\infty}(Q^+_-)$ and $\nabla u\in L_\infty(Q^+_-)$. The last fact has been proven in \cite{SS1}. Hence,
$${\rm div H}\in L_{2,\infty}(Q^+_-).$$
By the properties of the kernels $G$ and $K$, such a solution $\widetilde{u}$ is bounded and satisfies the energy identity
$$\int\limits_{\mathbb R^2_+}|\widetilde{u}(x,t)|^2dx +2\int\limits_A^t\int\limits_{\mathbb R^2_+}|\nabla \widetilde{u}(x,\tau)|^2dx d\tau=\int\limits_{\mathbb R^2_+}|{u}(x,A)|^2dx+$$
$$+2\int\limits_A^t\int\limits_{\mathbb R^2_+}{\rm div H}(x,\tau)\cdot \widetilde{u}(x,\tau)dx d\tau$$
for all $A\leq t\leq 0$. In addition,
 we can state that
 for any $\delta >0$,
 \begin{equation}\label{tildap2}
 \int\limits_{A+\delta}^0\int\limits_{\mathbb R^2_+}|\nabla \widetilde{p}^2|^2dx dt<C( \delta,A)<\infty.
 \end{equation}

Our aim is to show that $u=\widetilde{u}$ in $\mathbb R^2_+\times ]A,0[$. 
It is easy to see that,
for any $R>0$,
$$\|v(\cdot, t)\|_{2,B_+(R)}\to 0$$
as $t\to A$, where $v=u-\widetilde{u}$. This follows from the facts that $u$ is continuous on the completion of the set $Q_+(R)$ for any $R>0$, see details in \cite{SS1}, and that $\widetilde{u}\in C([A,0];L_2(\mathbb R^2_+))$.

The latter property allows us to
show that $v$ satisfies that the identity
$$\int\limits^0_A\int\limits_{\mathbb R^2_+}(v\cdot \partial_t\varphi+v\cdot \triangle\varphi)dxdt=0$$
for any $\varphi\in C^\infty_0(Q_-)$ such that  $\varphi(x_1,0,t)=0$  for any $x_1\in\mathbb R$ and any $t\in ]-\infty,0[$ and ${\rm div}\, \varphi=0$ in $Q^+_-$.
If we extend $v$ by zero for $t<A$, this field will be bounded ancient solution to the Stokes system and therefore has the form $v=(v_1(x_2,t),0)$, see \cite{JSS1} and \cite{JSS2}. The  gradient of the corresponding pressure $p^2-\widetilde{p}^2$ depends only on $t$. However, by (\ref{pressure2property})  and by (\ref{tildap2}), this gradient must be zero.
And the Liouville theorem for the heat equation in the half-space implies that $v=0$. 

Now, since $u=\widetilde{u}$, the energy identity implies
$$\int\limits_{\mathbb R^2_+}|{u}(x,0)|^2dx +2\int\limits_A^0\int\limits_{\mathbb R^2_+}|\nabla {u}(x,\tau)|^2dx d\tau=\int\limits_{\mathbb R^2_+}|{u}(x,A)|^2dx$$
for any $A<0$. This completes the proof of the lemma. $\Box$

\begin{remark}\label{uniformboundedness}
In fact, we have proven that
$$\int\limits^0_{A}\int\limits_{\mathbb R^2_+}|\nabla p^2|^2dx dt\leq c<\infty$$
for any $A<0$.
\end{remark}

Given a tensor-valued function $F\in C^\infty_0(Q_-^+)$,
let us consider the following initial boundary value problem:
\begin{equation}\label{perturbproblem}
    \partial_tv+u\cdot \nabla v+\triangle v+\nabla q={\rm div}\, F,\qquad {\rm div}\,v=0
\end{equation}
in $Q_+=\mathbb R^2_+\times ]-\infty,0[$,
\begin{equation}\label{boundarypertproblem}
  v(x_1,0,t)=0
\end{equation}
for any $x_1\in \mathbb R$ and $t\leq 0$, and
\begin{equation}\label{initialpertproblem}
    v(x,0)=0
\end{equation}
for $x\in \mathbb R^2_+$. Here, vector-valued field $v$ and scalar function $q$ are unknown.

Why we consider this system?
At least formally, we have the following identity
$$\int\limits_{Q_-^+}u\cdot {\rm div}\,{F}dx dt=$$
$$=\int\limits_{Q_-^+}u\cdot \Big( \partial_t{v}+u\cdot \nabla {v}+\triangle {v}+\nabla {q}\Big)dx dt =$$
$$=\int\limits_{Q_-^+}u\cdot \Big( \partial_t{v}+u\cdot \nabla {v}+\triangle {v}\Big)dx dt =$$
$$=\int\limits_{Q_-^+}\Big(-\partial_tu-{\rm div}\,u\otimes u+\triangle u\Big)\cdot {v}dx dt=$$
$$=\int\limits_{Q_-^+}\Big(-\partial_tu-{\rm div}\,u\otimes u+\triangle u-\nabla p\Big)\cdot {v}dx dt=0.$$
This would imply that $u$ is the function of $t$ only and thus, since $u$ is a mild bounded ancient solution, $u$ must be identically zero.

\pagebreak
\setcounter{equation}{0}
\section{Properties of Solutions to Dual Problem}
\begin{pro}\label{2p1} There exists a unique solution $v$ to (\ref{perturbproblem}), (\ref{boundarypertproblem}), and (\ref{initialpertproblem}) with the following properties:
$$v\in L_{2,\infty}(Q_-^+),\qquad \nabla v\in L_2(Q_-^+),$$
and, for all $T<0$,
$$\partial_tp, \nabla^2v, \nabla q\in L_2(\mathbb R^2_+\times ]T,0[).$$
\end{pro}
\textsc{Proof} 

First of all, there exists a unique energy solution. This follows from the identity
$$\int\limits_{Q_-^+}(u\cdot\nabla v)\cdot vdx dt=0$$
and from the inequality
$$\Big|-\int\limits_{Q_-^+}{\rm div}\,F\cdot vdx dt\Big|=\Big|\int\limits_{Q_-^+}F:\nabla v dx dt\Big|\leq \Big(\int\limits_{Q_-^+}|F|^2dx dt\Big)^\frac 12\Big(\int\limits_{Q_-^+}|\nabla v|^2dx dt\Big)^\frac 12$$
So, we can state that
\begin{equation}\label{energyestimate}
    v\in L_{2,\infty}(Q_-^+),\qquad \nabla v\in L_2(Q_-^+).
\end{equation}
The latter means that $u\cdot \nabla v\in L_2(Q_-^+)$. So, statements of Proposition \ref{2p1} follows from the theory for Stokes system.
\begin{flushright}
  $\Box$
\end{flushright}

\pagebreak

\setcounter{equation}{0}
\section{Main Formula, Integration by Parts}

For  smooth function $\psi\in C^\infty_0(\mathbb R^2\times \mathbb R)$, we have
$$\int\limits_{Q_-^+}u\cdot \psi{\rm div}\,{F}dx dt=$$
$$=\int\limits_{Q_-^+}u\cdot \psi\Big(\partial_t{v}+u\cdot\nabla {v}+\triangle {v}+\nabla {q}\Big)dx dt=$$$$=\int\limits_{Q_-^+}\Big(-u\cdot {v}\partial_t\psi-u\cdot {v}u\cdot\nabla\psi-u_i{v}_{i,j}\psi_{,j}+u_{i,j}
{v}_i\psi_{,j}-{q}u\cdot\nabla\psi\Big)dx dt-$$
$$-{v}\psi\cdot\Big(\partial_tu+u\cdot\nabla u-\triangle u\Big)dx dt=$$
$$=\int\limits_{Q_-^+}\Big(-u\cdot {v}\partial_t\psi-u\cdot {v}u\cdot\nabla\psi-2u_i{v}_{i,j}\psi_{,j}+
(u_{i,j}{v}_i+u_i{v}_{i,j})\psi_{,j}-{q}u\cdot\nabla\psi\Big)dx dt+$$$$+\int\limits_{Q_-^+}{v}\psi\cdot \nabla pdx dt=$$
$$=-\int\limits_{Q_-^+}\Big(u\cdot {v}\partial_t\psi+u\cdot {v}u\cdot\nabla\psi+2u_i{v}_{i,j}\psi_{,j}+u\cdot {v} \triangle\psi
+({q}u+p{v})\cdot\nabla\psi\Big)dx dt.$$

We pick $\psi(x,t)=\chi(t)\varphi(x)$. Using simple arguments and smoothness of $u$ and $v$, we can get rid of $\chi$ and have
$$J(T)=\int\limits_T^0\int\limits_{\mathbb R^2_+}u\cdot \varphi{\rm div}\,{F}dx dt=-\int\limits_{\mathbb R^2_+}\varphi(x)u(x,T)\cdot {v}(x,T)dx+$$$$+\int\limits_T^0\int\limits_{\mathbb R^2_+}\Big(u\cdot {v}u\cdot\nabla\varphi+2u_i{v}_{i,j}\varphi_{,j}+u\cdot {v} \triangle\varphi
+({q}u+p{v})\cdot\nabla\varphi\Big)dx dt.$$

Fix a cut-off function
$\varphi(x)=\xi(x/R)$, where $\xi\in C^\infty_0(\mathbb R^3)$ with the following properties: $0\leq \xi\leq 1$, $\xi(x)=1$ if $|x|\leq 1$, and $\xi(x)=0$ if $|x|\geq 2$.
Our aim is to show that
$$J_R=\int\limits_T^0\int\limits_{\mathbb R^2_+}\Big(u\cdot \widetilde{v}u\cdot\nabla\varphi+2u_i\widetilde{v}_{i,j}\varphi_{,j}+u\cdot \widetilde{v} \triangle\varphi
+(\widetilde{q}u+p\widetilde{v})\cdot\nabla\varphi\Big)dx dt$$
tends to zero if $R\to \infty$.

We start with
$$\Big|\int\limits_T^0\int\limits_{\mathbb R^2_+}2u_i{v}_{i,j}\varphi_{,j}dx dt\Big|\leq $$$$\leq \frac cR\Big(\int\limits_T^0\int\limits_{B_+(2R)}|u|^2dx dt\Big)^\frac 12\Big(\int\limits_T^0\int\limits_{B_+(2R)\setminus B_+(R)}|\nabla {v}|^2dx dt\Big)^\frac 12\leq $$$$\leq c\sqrt {-T}\Big(\int\limits_T^0\int\limits_{B_+(2R)\setminus B_+(R)}|\nabla {v}|^2dx dt\Big)^\frac 12  \to 0$$
as $R\to \infty$.

Next, since $|u|\leq 1$, we have
$$\Big|\int\limits_T^0\int\limits_{\mathbb R^2_+}u\cdot {v} \triangle\varphi dx dt \Big|\leq \frac c{R^2}\Big(\int\limits_T^0\int\limits_{B_+(2R)}|u|^2dx dt\Big)^\frac 12\Big(\int\limits_T^0\int\limits_{\mathbb R^2_+}| {v}|^2dx dt\Big)^\frac 12\leq $$$$\leq c  \frac {-T}{R}\|{v}\|_{2,\infty,Q_-^+}\to 0$$
as $R\to \infty$.

The third term is estimated as follows (by boundedness of $u$):
$$\Big|\int\limits_T^0\int\limits_{\mathbb R^2_+}u\cdot {v}u\cdot\nabla\varphi dx dt\Big|\leq$$$$\leq  \frac c{R}\Big(\int\limits_T^0\int\limits_{B_+(2R)}|u|^4dx dt\Big)^\frac 12\Big(\int\limits_T^0\int\limits_{B_+(2R)\setminus B_+(R)}| {v}|^2dx dt\Big)^\frac 12 \leq$$
$$\leq  \frac c{R}\Big(\int\limits_T^0\int\limits_{B_+(2R)}|u|^2dx dt\Big)^\frac 12\Big(\int\limits_T^0\int\limits_{B_+(2R)\setminus B_+(R)}| {v}|^2dx dt\Big)^\frac 12 \leq$$$$\leq  \frac c{R}\Big(\int\limits_T^0\int\limits_{\mathbb R^2_+}|u|^2dx dt\Big)^\frac 12\Big(\int\limits_T^0\int\limits_{\mathbb R^2_+}| {v}|^2dx dt\Big)^\frac 12 \leq$$$$\leq c\frac {-T}R\|u\|_{2,\infty,Q^+_-}\|u\|_{2,\infty,Q^+_-}\to 0$$
as $R\to \infty$.

The first term containing the  pressure is estimated as follows. We have
$$\int\limits_T^0\int\limits_{\mathbb R^2_+}p{v}\cdot\nabla\varphi dx dt=
\int\limits_T^0\int\limits_{\mathbb R^2_+}p_R{v}\cdot\nabla\varphi dx dt,$$
where
$$p_R=p^1_R+p^2_R$$
with $p^1_R=p^1-[p^1]_{B_+(2R)}$ and $p^2_R=p^2-[p^2]_{B_+(2R)}$. By the assumptions,
 after even extension, the function $p^1$ belongs to $L_\infty(-\infty, 0;BMO)$ and thus
$$\frac 1{R^2}\int\limits_{B_+(2R)}|p^1_R(x,t)|^2 dx\leq c$$
for all $t\leq 0$.  As to $p^2_R$, we use Poincar\'{e} inequality
$$\frac 1{R^2}\int\limits_{B_+(2R)}|p^2_R(x,t)|^2 dx\leq \int\limits_{B_+(2R)}|\nabla p^2(x,t)|^2dx\leq \int\limits_{\mathbb R^2_+}|\nabla p^2(x,t)|^2dx.$$

So, by Lemma \ref{dissipationisbounded} and by the Lebesgue theorem about dominated convergence,
$$\Big|\int\limits_T^0\int\limits_{\mathbb R^2_+}p{v}\cdot\nabla\varphi dx dt\Big|\leq $$$$
\leq c\int\limits_T^0d\tau\Big(\int\limits_{B_+(2R)\setminus B_+(R)}|v(x,\tau)|^2dx\Big)^\frac 12+$$$$
+\frac cR\Big(R^2\int\limits^0_T\int\limits_{\mathbb R^2_+}|\nabla p^2|^2dx dt\Big)^\frac 12\Big(\int\limits^0_T\int\limits_{B_+(2R)\setminus B_+(R)}|v|^2dx dt\Big)^\frac 12
\to 0$$
as $R\to \infty$.

The last term is  treated with the help of  Poincar\'{e} inequality in the same way as $p^2_R$. Indeed,
$$\Big|\int\limits_T^0\int\limits_{\mathbb R^2_+}qu\cdot\nabla\varphi dx dt\Big|
=\Big|\int\limits_T^0\int\limits_{\mathbb R^2_+}(q-[q]_{B_+(2R)})u\cdot\nabla\varphi dx dt\Big|\leq$$
$$\leq \frac cR\Big(R^2\int\limits^0_T\int\limits_{B_+(2R)}|\nabla q|^2dx dt\Big)^\frac 12\Big(\int\limits^0_T\int\limits_{B_+(2R)\setminus B_+(R)}|u|^2dx dt\Big)^\frac 12.$$
The right hand side of the latter inequality tends to zero as $R\to\infty$ by the assumption that $u\in L_{2,\infty}(Q_-^+)$.

So, finally, we have
$$\int\limits^0_T\int\limits_{\mathbb R^2_+} u\cdot {\rm div}\,F dx dt=
-\lim\limits_{R\to\infty}\int\limits_{\mathbb R^2_+}\varphi(x)u(x,T)\cdot {v}(x,T)dx=$$
$$=-\int\limits_{\mathbb R^2_+}u(x,T)\cdot {v}(x,T)dx$$

Now, our aim to see what happens if $T\to-\infty$.


\pagebreak
\setcounter{equation}{0}
\section{$t\to-\infty$}

We shall show that
\begin{equation}\label{l2decay}
    \|v(\cdot,t)\|_{2,\mathbb R^2_+}\to 0.
\end{equation}
as $t\to-\infty$.

Indeed, we also know
\begin{equation}\label{dissipto0}
  \int\limits_{-\infty}^t\int\limits_{\mathbb R^2_+}|\nabla v|^2dx d\tau\to0
\end{equation}
as $t\to-\infty$.
By Ladyzhenskaya's inequality,
$$v\in L_4(Q_-^+)$$
and thus
\begin{equation}\label{L}
  \int\limits_{-\infty}^t\int\limits_{\mathbb R^2_+}| v|^4dx d\tau\to0
\end{equation}
as $t\to-\infty$.

Now, for sufficiently large $-t_0$, we have
$$v=v^1+v^2,$$
where
$$\partial_tv^1 +\triangle v^1+\nabla q^1=0,\qquad {\rm div}\,v^1=0$$
in $R^2_+\times ]-\infty,t_0[$,
$$v^1(x_1,0,t)=0$$
for any $x_1\in\mathbb R$ and for any $t\leq t_0$, and
$$v^1(x,t_0)=v(x,t_0)$$
for any $x\in \mathbb R^2_+$.

As to $v^2$, it satisfies
$$\partial_tv^2 +\triangle v^2+\nabla q^2=-{\rm div}\, v\otimes u,\qquad {\rm div}\,v^2=0$$
in $R^2_+\times ]-\infty,t_0[$,
$$v^2(x_1,0,t)=0$$
for any $x_1\in\mathbb R$ and for any $t\leq t_0$, and
$$v^1(x,t_0)=0$$
for any $x\in \mathbb R^2_+$.

Then, it is well known that
\begin{equation}\label{l2decay1}    \|v^1(\cdot,t)\|_{2,\mathbb R^2_+}\to 0\end{equation}
as $t\to-\infty$.
On the other hand, by the energy inequality,
$$ \frac 12 \|v^2(\cdot,t)\|_{2,\mathbb R^2_+}^2+\int\limits^{t_0}_t\int\limits_{\mathbb R^2_+}|\nabla v^2|^2dx d\tau=$$$$= \int\limits_t^{t_0}\int\limits_{\mathbb R^2_+}v^2\cdot{\rm div} \,v\otimes u dx d\tau=
- \int\limits_t^{t_0}\int\limits_{\mathbb R^2_+} \,v\otimes u:\nabla v^2 dx d\tau\leq$$
$$\leq \Big(\int\limits^{t_0}_t\int\limits_{\mathbb R^2_+}|\nabla v^2|^2dx d\tau\Big)^\frac 12\Big(\int\limits^{t_0}_t\int\limits_{\mathbb R^2_+}|u|^4dx d\tau\Big)^\frac 14\Big(\int\limits^{t_0}_t\int\limits_{\mathbb R^2_+}| v|^4dx d\tau\Big)^\frac 14$$
for $t<t_0$.
For the same reason as for $v$, we have
\begin{equation}\label{essential}
  \int\limits^{0}_{\infty}\int\limits_{\mathbb R^2_+}|u|^4dx d\tau\leq c
\end{equation}

and thus, by the Cauchy inequality,
\begin{equation}\label{l2decay2}\|v^2(\cdot,t)\|_{2,\mathbb R^2_+}^2\leq c\Big(\int\limits^{t_0}_t\int\limits_{\mathbb R^2_+}| v|^4dx d\tau\Big)^\frac 12\leq c\Big(\int\limits^{t_0}_{-\infty}\int\limits_{\mathbb R^2_+}| v|^4dx d\tau\Big)^\frac 12\end{equation}
for all $t<t_0$.

It is not so difficult to deduce (\ref{l2decay}) from (\ref{dissipto0}), (\ref{l2decay1}), and (\ref{l2decay2}).

The only assumption we really need is (\ref{essential}) and it is true if $u\in L_{2,\infty}(Q^+_-)$ and $\nabla u\in L_2(Q^+_-)$. The latter follows from Ladyzhenskaya's inequality.

\pagebreak

\end{document}